\documentclass{article}
\usepackage{amsfonts}
\usepackage{amsmath}
\usepackage{amsthm}

\theoremstyle{plain}
\newtheorem{theorem}{Theorem}
\theoremstyle{definition}

\newtheorem{definition}[theorem]{Definition}
\newtheorem{example}[theorem]{Example}
\newtheorem{remark}[theorem]{Remark}
\newenvironment{demo}[1][Proof]{\noindent\textbf{#1.} }{\hfill$\Box$}
\numberwithin{theorem}{section}
\numberwithin{equation}{section}

\input{tcilatex}

\begin{document}

\title{Jet geometrical extension of the KCC-invariants}
\author{Vladimir Balan and Mircea Neagu \and {\small June 2009; Last revised
December 2009 (an added bibliographical item)}}
\date{}
\maketitle

\begin{abstract}
In this paper we construct the jet geometrical extensions of the
KCC-invariants, which characterize a given second-order system of
differential equations on the 1-jet space $J^{1}(\mathbb{R} ,M)$. A
generalized theorem of characterization of our jet geometrical
KCC-invariants is also presented.
\end{abstract}

\textbf{Mathematics Subject Classification (2000):} 58B20, 37C10, 53C43.

\textbf{Key words and phrases:} 1-jet spaces, temporal and spatial
semisprays, nonlinear connections, SODEs, jet $h-$KCC-invariants.

\section{Geometrical objects on 1-jet spaces}

\hspace{5mm}We remind first several differential geometrical properties of
the 1-jet spaces. The 1-jet bundle%
\begin{equation*}
\xi =(J^{1}(\mathbb{R},M),\pi _{1},\mathbb{R}\times M)
\end{equation*}%
is a vector bundle over the product manifold $\mathbb{R}\times M$, having
the fibre of type $\mathbb{R}^{n}$, where $n$ is the dimension of the 
\textit{spatial} manifold $M$. If the spatial manifold $M$ has the local
coordinates $(x^{i})_{i=\overline{1,n}}$, then we shall denote the local
coordinates of the 1-jet total space $J^{1}(\mathbb{R},M)$ by $%
(t,x^{i},x_{1}^{i})$; these transform by the rules \cite{8}%
\begin{equation}
\left\{ 
\begin{array}{l}
\widetilde{t}=\widetilde{t}(t)\medskip \\ 
\widetilde{x}^{i}=\widetilde{x}^{i}(x^{j})\medskip \\ 
\widetilde{x}_{1}^{i}=\dfrac{\partial \widetilde{x}^{i}}{\partial x^{j}}%
\dfrac{dt}{d\widetilde{t}}\cdot x_{1}^{j}.%
\end{array}%
\right.  \label{rgg}
\end{equation}%
In the geometrical study of the 1-jet bundle, a central role is played by
the \textit{distinguished tensors} ($d-$tensors).

\begin{definition}
A geometrical object $D=\left( D_{1k(1)(l)...}^{1i(j)(1)...}\right) $ on the
1-jet vector bundle, whose local components transform by the rules%
\begin{equation}
D_{1k(1)(u)...}^{1i(j)(1)...}=\widetilde{D}_{1r(1)(s)...}^{1p(m)(1)...}\frac{%
dt}{d\widetilde{t}}\frac{\partial x^{i}}{\partial \widetilde{x}^{p}}\left( 
\frac{\partial x^{j}}{\partial \widetilde{x}^{m}}\frac{d\widetilde{t}}{dt}%
\right) \frac{d\widetilde{t}}{dt}\frac{\partial \widetilde{x}^{r}}{\partial
x^{k}}\left( \frac{\partial \widetilde{x}^{s}}{\partial x^{u}}\frac{dt}{d%
\widetilde{t}}\right) ...,  \label{tr-rules-$d-$tensors}
\end{equation}%
is called a \textit{$d-$tensor field}.
\end{definition}

\begin{remark}
The use of parentheses for certain indices of the local components%
\begin{equation*}
D_{1k(1)(l)...}^{1i(j)(1)...}
\end{equation*}%
of the distinguished tensor field $D$ on the 1-jet space is motivated by the
fact that the pair of indices $"$ $_{(1)}^{(j)}$ $"$ or $"$ $_{(l)}^{(1)}$ $%
" $ behaves like a single index.
\end{remark}

\begin{example}
\label{Liouville} The geometrical object%
\begin{equation*}
\mathbf{C}=\mathbf{C}_{(1)}^{(i)}\dfrac{\partial }{\partial x_{1}^{i}},
\end{equation*}%
where $\mathbf{C}_{(1)}^{(i)}=x_{1}^{i}$, represents a $d-$tensor field on
the 1-jet space; this is called the \textit{canonical Liouville $d-$tensor
field} of the 1-jet bundle and is a global geometrical object.
\end{example}

\begin{example}
\label{normal}Let $h=(h_{11}(t))$ be a Riemannian metric on the relativistic
time axis $\mathbb{R}$. The geometrical object 
\begin{equation*}
\mathbf{J}_{h}=J_{(1)1j}^{(i)}\dfrac{\partial }{\partial x_{1}^{i}}\otimes
dt\otimes dx^{j},
\end{equation*}%
where $J_{(1)1j}^{(i)}=h_{11}\delta _{j}^{i}$ is a $d-$tensor field on $%
J^{1}(\mathbb{R},M)$, which is called the $h$\textit{-normalization $d-$%
tensor field} of the 1-jet space and is a global geometrical object.
\end{example}

In the Riemann-Lagrange differential geometry of the 1-jet spaces developed
in \cite{7}, \cite{8} important r\^oles are also played by geometrical
objects as the \textit{temporal} or \textit{spatial semisprays}, together
with the \textit{jet nonlinear connections}.

\begin{definition}
A set of local functions $H=\left( H_{(1)1}^{(j)}\right) $ on $J^{1}(\mathbb{%
R},M),$ which transform by the rules%
\begin{equation}
2\widetilde{H}_{(1)1}^{(k)}=2H_{(1)1}^{(j)}\left( \frac{dt}{d\widetilde{t}}%
\right) ^{2}\frac{\partial \widetilde{x}^{k}}{\partial x^{j}}-\frac{dt}{d%
\widetilde{t}}\frac{\partial \widetilde{x}_{1}^{k}}{\partial t},
\label{tr-rules-t-s}
\end{equation}%
is called a \textit{temporal semispray} on $J^{1}(\mathbb{R} ,M)$.
\end{definition}

\begin{example}
\label{H0} Let us consider a Riemannian metric $h=(h_{11}(t))$ on the
temporal manifold $\mathbb{R} $ and let%
\begin{equation*}
H_{11}^{1}=\frac{h^{11}}{2}\frac{dh_{11}}{dt},
\end{equation*}%
where $h^{11}=1/h_{11}$, be its Christoffel symbol. Taking into account that
we have the transformation rule%
\begin{equation}
\widetilde{H}_{11}^{1}=H_{11}^{1}\frac{dt}{d\widetilde{t}}+\frac{d\widetilde{%
t}}{dt}\frac{d^{2}t}{d\widetilde{t}^{2}},  \label{t-Cris-symb}
\end{equation}%
we deduce that the local components%
\begin{equation*}
\mathring{H}_{(1)1}^{(j)}=-\frac{1}{2}H_{11}^{1}x_{1}^{j}
\end{equation*}%
define a temporal semispray $\mathring{H}=\left( \mathring{H}%
_{(1)1}^{(j)}\right) $ on $J^{1}(\mathbb{R} ,M)$. This is called the \textit{%
canonical temporal semispray associated to the temporal metric} $h(t)$.
\end{example}

\begin{definition}
A set of local functions $G=\left( G_{(1)1}^{(j)}\right) ,$ which transform
by the rules%
\begin{equation}
2\widetilde{G}_{(1)1}^{(k)}=2G_{(1)1}^{(j)}\left( \frac{dt}{d\widetilde{t}}%
\right) ^{2}\frac{\partial \widetilde{x}^{k}}{\partial x^{j}}-\frac{\partial
x^{m}}{\partial \widetilde{x}^{j}}\frac{\partial \widetilde{x}_{1}^{k}}{%
\partial x^{m}}\widetilde{x}_{1}^{j},  \label{tr-rules-s-s}
\end{equation}%
is called a \textit{spatial semispray} on $J^{1}(\mathbb{R} ,M)$.
\end{definition}

\begin{example}
\label{G0} Let $\varphi =(\varphi _{ij}(x))$ be a Riemannian metric on the
spatial manifold $M$ and let us consider%
\begin{equation*}
\gamma _{jk}^{i}=\frac{\varphi ^{im}}{2}\left( \frac{\partial \varphi _{jm}}{%
\partial x^{k}}+\frac{\partial \varphi _{km}}{\partial x^{j}}-\frac{\partial
\varphi _{jk}}{\partial x^{m}}\right)
\end{equation*}%
its Christoffel symbols. Taking into account that we have the transformation
rules%
\begin{equation}
\widetilde{\gamma }_{qr}^{p}=\gamma _{jk}^{i}\frac{\partial \widetilde{x}^{p}%
}{\partial x^{i}}\frac{\partial x^{j}}{\partial \widetilde{x}^{q}}\frac{%
\partial x^{k}}{\partial \widetilde{x}^{r}}+\frac{\partial \widetilde{x}^{p}%
}{\partial x^{l}}\frac{\partial ^{2}x^{l}}{\partial \widetilde{x}%
^{q}\partial \widetilde{x}^{r}},  \label{s-Cris-symb}
\end{equation}%
we deduce that the local components%
\begin{equation*}
\mathring{G}_{(1)1}^{(j)}=\frac{1}{2}\gamma _{kl}^{j}x_{1}^{k}x_{1}^{l}
\end{equation*}%
define a spatial semispray $\mathring{G}=\left( \mathring{G}%
_{(1)1}^{(j)}\right) $ on $J^{1}(\mathbb{R} ,M)$. This is called the \textit{%
canonical spatial semispray associated to the spatial metric} $\varphi (x)$.
\end{example}

\begin{definition}
A set of local functions $\Gamma =\left(
M_{(1)1}^{(j)},N_{(1)i}^{(j)}\right) $ on $J^{1}(\mathbb{R} ,M),$ which
transform by the rules%
\begin{equation}
\widetilde{M}_{(1)1}^{(k)}=M_{(1)1}^{(j)}\left( \frac{dt}{d\widetilde{t}}%
\right) ^{2}\frac{\partial \widetilde{x}^{k}}{\partial x^{j}}-\frac{dt}{d%
\widetilde{t}}\frac{\partial \widetilde{x}_{1}^{k}}{\partial t}
\label{tr-rules-t-nlc}
\end{equation}%
and%
\begin{equation}
\widetilde{N}_{(1)l}^{(k)}=N_{(1)i}^{(j)}\frac{dt}{d\widetilde{t}}\frac{%
\partial x^{i}}{\partial \widetilde{x}^{l}}\frac{\partial \widetilde{x}^{k}}{%
\partial x^{j}}-\frac{\partial x^{m}}{\partial \widetilde{x}^{l}}\frac{%
\partial \widetilde{x}_{1}^{k}}{\partial x^{m}},  \label{tr-rules-s-nlc}
\end{equation}%
is called a \textit{nonlinear connection} on the 1-jet space $J^{1}(\mathbb{R%
},M)$.
\end{definition}

\begin{example}
Let us consider that $(\mathbb{R} ,h_{11}(t))$ and $(M,\varphi _{ij}(x))$
are Rie\-ma\-nni\-an manifolds having the Christoffel symbols $H_{11}^{1}(t)$
and $\gamma _{jk}^{i}(x)$. Then, using the transformation rules (\ref{rgg}),
(\ref{t-Cris-symb}) and (\ref{s-Cris-symb}), we deduce that the set of local
functions%
\begin{equation*}
\mathring{\Gamma}=\left( \mathring{M}_{(1)1}^{(j)},\mathring{N}%
_{(1)i}^{(j)}\right) ,
\end{equation*}%
where%
\begin{equation*}
\mathring{M}_{(1)1}^{(j)}=-H_{11}^{1}x_{1}^{j}\text{ \ \ and \ \ }\mathring{N%
}_{(1)i}^{(j)}=\gamma _{im}^{j}x_{1}^{m},
\end{equation*}%
represents a nonlinear connection on the 1-jet space $J^{1}(\mathbb{R} ,M)$.
This jet nonlinear connection is called the \textit{canonical nonlinear
connection attached to the pair of Riemannian metrics} $(h(t),\varphi (x))$.
\end{example}

In the sequel, let us study the geometrical relations between \textit{%
temporal }or\textit{\ spatial semisprays} and \textit{nonlinear connections}
on the 1-jet space $J^{1}(\mathbb{R} ,M)$. In this direction, using the
local transformation laws (\ref{tr-rules-t-s}), (\ref{tr-rules-t-nlc}) and (%
\ref{rgg}), respectively the transformation laws (\ref{tr-rules-s-s}), (\ref%
{tr-rules-s-nlc}) and (\ref{rgg}), by direct local computation, we find the
following geometrical results:

\begin{theorem}
a) The \textit{temporal semisprays} $H=(H_{(1)1}^{(j)})$ and the sets of 
\textit{temporal components of nonlinear connections }$\Gamma _{\text{%
temporal}}=(M_{(1)1}^{(j)})$ are in one-to-one correspondence on the 1-jet
space $J^{1}(\mathbb{R},M)$, via: 
\begin{equation*}
M_{(1)1}^{(j)}=2H_{(1)1}^{(j)},\qquad H_{(1)1}^{(j)}=\frac{1}{2}%
M_{(1)1}^{(j)}.
\end{equation*}

b) The \textit{spatial semisprays} $G=(G_{(1)1}^{(j)})$ and the sets of 
\textit{spatial components of nonlinear connections }$\Gamma _{\text{spatial}%
}=(N_{(1)k}^{(j)})$ are connected on the 1-jet space $J^{1}(\mathbb{R},M)$,
via the relations: 
\begin{equation*}
N_{(1)k}^{(j)}=\frac{\partial G_{(1)1}^{(j)}}{\partial x_{1}^{k}},\qquad
G_{(1)1}^{(j)}=\frac{1}{2}N_{(1)m}^{(j)}x_{1}^{m}.
\end{equation*}
\end{theorem}

\section{Jet geometrical KCC-theory}

\hspace{5mm}In this Section we generalize on the 1-jet space $J^{1}(\mathbb{R%
},M)$ the basics of the KCC-theory (\cite{1}, \cite{2}, \cite{3}, \cite{9}).
In this respect, let us consider on $J^{1}(\mathbb{R} ,M)$ a second-order
system of differential equations of local form%
\begin{equation}
\frac{d^{2}x^{i}}{dt^{2}}+F_{(1)1}^{(i)}(t,x^{k},x_{1}^{k})=0,\text{ \ \ }i=%
\overline{1,n},  \label{SODE}
\end{equation}%
where $x_{1}^{k}=dx^{k}/dt$ and the local components $%
F_{(1)1}^{(i)}(t,x^{k},x_{1}^{k})$ transform under a change of coordinates (%
\ref{rgg}) by the rules%
\begin{equation}
\widetilde{F}_{(1)1}^{(r)}=F_{(1)1}^{(j)}\left( \frac{dt}{d\widetilde{t}}%
\right) ^{2}\frac{\partial \widetilde{x}^{r}}{\partial x^{j}}-\frac{dt}{d%
\widetilde{t}}\frac{\partial \widetilde{x}_{1}^{r}}{\partial t}-\frac{%
\partial x^{m}}{\partial \widetilde{x}^{j}}\frac{\partial \widetilde{x}%
_{1}^{r}}{\partial x^{m}}\widetilde{x}_{1}^{j}.  \label{transformations-F}
\end{equation}

\begin{remark}
The second-order system of differential equations (\ref{SODE}) is invariant
under a change of coordinates (\ref{rgg}).
\end{remark}

Using a temporal Riemannian metric $h_{11}(t)$ on $\mathbb{R} $ and taking
into account the transformation rules (\ref{tr-rules-t-s}) and (\ref%
{tr-rules-s-s}), we can rewrite the SODEs (\ref{SODE}) in the following form:%
\begin{equation*}
\frac{d^{2}x^{i}}{dt^{2}}%
-H_{11}^{1}x_{1}^{i}+2G_{(1)1}^{(i)}(t,x^{k},x_{1}^{k})=0,\text{ \ \ }i=%
\overline{1,n},
\end{equation*}%
where%
\begin{equation*}
G_{(1)1}^{(i)}=\frac{1}{2}F_{(1)1}^{(i)}+\frac{1}{2}H_{11}^{1}x_{1}^{i}
\end{equation*}%
are the components of a spatial semispray on $J^{1}(\mathbb{R} ,M)$.
Moreover, the coefficients of the spatial semispray $G_{(1)1}^{(i)}$ produce
the spatial components $N_{(1)j}^{(i)}$ of a nonlinear connection $\Gamma $
on the 1-jet space $J^{1}(\mathbb{R} ,M)$, by putting%
\begin{equation*}
N_{(1)j}^{(i)}=\frac{\partial G_{(1)1}^{(i)}}{\partial x_{1}^{j}}=\frac{1}{2}%
\frac{\partial F_{(1)1}^{(i)}}{\partial x_{1}^{j}}+\frac{1}{2}%
H_{11}^{1}\delta _{j}^{i}.
\end{equation*}

In order to find the basic jet differential geometrical invariants of the
system (\ref{SODE}) (see Kosambi \cite{6}, Cartan \cite{4} and Chern \cite{5}%
) under the jet coordinate transformations (\ref{rgg}), we define the $h-$%
\textit{KCC-covariant derivative of a $d-$tensor of kind }$%
T_{(1)}^{(i)}(t,x^{k},x_{1}^{k})$ on the 1-jet space $J^{1}(\mathbb{R} ,M)$
via%
\begin{eqnarray*}
\frac{\overset{h}{D}T_{(1)}^{(i)}}{dt} &=&\frac{dT_{(1)}^{(i)}}{dt}%
+N_{(1)r}^{(i)}T_{(1)}^{(r)}-H_{11}^{1}T_{(1)}^{(i)}= \\
&=&\frac{dT_{(1)}^{(i)}}{dt}+\frac{1}{2}\frac{\partial F_{(1)1}^{(i)}}{%
\partial x_{1}^{r}}T_{(1)}^{(r)}-\frac{1}{2}H_{11}^{1}T_{(1)}^{(i)},
\end{eqnarray*}%
where the Einstein summation convention is used throughout.

\begin{remark}
The $h-$\textit{KCC-covariant derivative} components $\dfrac{\overset{h}{D}%
T_{(1)}^{(i)}}{dt}$ transform under a change of coordinates (\ref{rgg}) as a 
$d-$tensor of type $\mathcal{T}_{(1)1}^{(i)}.$
\end{remark}

In such a geometrical context, if we use the notation $x_{1}^{i}=dx^{i}/dt$,
then the system (\ref{SODE}) can be rewritten in the following distinguished
tensorial form:%
\begin{eqnarray*}
\frac{\overset{h}{D}x_{1}^{i}}{dt}
&=&-F_{(1)1}^{(i)}(t,x^{k},x_{1}^{k})+N_{(1)r}^{(i)}x_{1}^{r}-H_{11}^{1}x_{1}^{i}=
\\
&=&-F_{(1)1}^{(i)}+\frac{1}{2}\frac{\partial F_{(1)1}^{(i)}}{\partial
x_{1}^{r}}x_{1}^{r}-\frac{1}{2}H_{11}^{1}x_{1}^{i},
\end{eqnarray*}

\begin{definition}
The distinguished tensor%
\begin{equation*}
\overset{h}{\varepsilon }\text{ }\!\!_{(1)1}^{(i)}=-F_{(1)1}^{(i)}+\frac{1}{2%
}\frac{\partial F_{(1)1}^{(i)}}{\partial x_{1}^{r}}x_{1}^{r}-\frac{1}{2}%
H_{11}^{1}x_{1}^{i}
\end{equation*}%
is called the \textit{first }$h-$\textit{KCC-invariant} on the 1-jet space $%
J^{1}(\mathbb{R} ,M)$ of the SODEs (\ref{SODE}), which is interpreted as an 
\textit{external force} \cite{1}, \cite{3}.
\end{definition}

\begin{example}
It can be easily seen that for the particular first order jet rheonomic
dynamical system%
\begin{equation}
\frac{dx^{i}}{dt}=X_{(1)}^{(i)}(t,x^{k})\Rightarrow \frac{d^{2}x^{i}}{dt^{2}}%
=\frac{\partial X_{(1)}^{(i)}}{\partial t}+\frac{\partial X_{(1)}^{(i)}}{%
\partial x^{m}}x_{1}^{m},  \label{Jet_DS}
\end{equation}%
where $X_{(1)}^{(i)}(t,x)$ is a given $d-$tensor on $J^{1}(\mathbb{R} ,M)$,
the first\textit{\ }$h-$KCC-invariant has the form%
\begin{equation*}
\overset{h}{\varepsilon }\text{ }\!\!_{(1)1}^{(i)}=\frac{\partial
X_{(1)}^{(i)}}{\partial t}+\frac{1}{2}\frac{\partial X_{(1)}^{(i)}}{\partial
x^{r}}x_{1}^{r}-\frac{1}{2}H_{11}^{1}x_{1}^{i}.
\end{equation*}
\end{example}

In the sequel, let us vary the trajectories $x^{i}(t)$ of the system (\ref%
{SODE}) by the nearby trajectories $(\overline{x}^{i}(t,s))_{s\in
(-\varepsilon ,\varepsilon )},$ where $\overline{x}^{i}(t,0)=x^{i}(t).$
Then, considering the \textit{variation $d-$tensor field}%
\begin{equation*}
\mathit{\ }\xi ^{i}(t)=\left. \dfrac{\partial \overline{x}^{i}}{\partial s}%
\right\vert _{s=0},
\end{equation*}%
we get the \textit{variational equations}%
\begin{equation}
\frac{d^{2}\xi ^{i}}{dt^{2}}+\frac{\partial F_{(1)1}^{(i)}}{\partial x^{j}}%
\xi ^{j}+\frac{\partial F_{(1)1}^{(i)}}{\partial x_{1}^{r}}\frac{d\xi ^{r}}{%
dt}=0.  \label{Var-Equations}
\end{equation}

In order to find other jet geometrical invariants for the system (\ref{SODE}%
), we also introduce the $h-$\textit{KCC-covariant derivative of a $d-$%
tensor of kind }$\xi ^{i}(t)$ on the 1-jet space $J^{1}(\mathbb{R} ,M)$ via 
\begin{equation*}
\frac{\overset{h}{D}\xi ^{i}}{dt}=\frac{d\xi ^{i}}{dt}+N_{(1)m}^{(i)}%
\xi^{m}= \frac{d\xi ^{i}}{dt}+\frac{1}{2}\frac{\partial F_{(1)1}^{(i)}}{%
\partial x_{1}^{m}} \xi ^{m}+\frac{1}{2}H_{11}^{1}\xi ^{i}.
\end{equation*}

\begin{remark}
The $h-$\textit{KCC-covariant derivative} components $\dfrac{\overset{h}{D}%
\xi ^{i}}{dt}$ transform under a change of coordinates (\ref{rgg}) as a $d-$%
tensor $T_{(1)}^{(i)}.$
\end{remark}

In this geometrical context, the variational equations (\ref{Var-Equations})
can be rewritten in the following distinguished tensorial form:%
\begin{equation*}
\frac{\overset{h}{D}}{dt}\left[ \frac{\overset{h}{D}\xi ^{i}}{dt}\right] =%
\overset{h}{P}\text{ \negthinspace \negthinspace }_{m11}^{i}\xi ^{m},
\end{equation*}%
where 
\begin{eqnarray*}
\overset{h}{P}\text{ \negthinspace \negthinspace }_{j11}^{i} &=&-\frac{%
\partial F_{(1)1}^{(i)}}{\partial x^{j}}+\frac{1}{2}\frac{\partial
^{2}F_{(1)1}^{(i)}}{\partial t\partial x_{1}^{j}}+\frac{1}{2}\frac{\partial
^{2}F_{(1)1}^{(i)}}{\partial x^{r}\partial x_{1}^{j}}x_{1}^{r}-\frac{1}{2}%
\frac{\partial ^{2}F_{(1)1}^{(i)}}{\partial x_{1}^{r}\partial x_{1}^{j}}%
F_{(1)1}^{(r)}+ \\
&&+\frac{1}{4}\frac{\partial F_{(1)1}^{(i)}}{\partial x_{1}^{r}}\frac{%
\partial F_{(1)1}^{(r)}}{\partial x_{1}^{j}}+\frac{1}{2}\frac{dH_{11}^{1}}{dt%
}\delta _{j}^{i}-\frac{1}{4}H_{11}^{1}H_{11}^{1}\delta _{j}^{i}.
\end{eqnarray*}

\begin{definition}
The $d-$tensor $\overset{h}{P}$ \negthinspace $_{j11}^{i}$ is called the 
\textit{second }$h-$\textit{KCC-invariant} on the 1-jet space $J^{1}(\mathbb{%
R},M)$ of the system (\ref{SODE}), or the \textit{jet }$h-$\textit{deviation
curvature $d-$tensor}.
\end{definition}

\begin{example}
If we consider the second-order system of differential equations of the 
\textit{harmonic curves associated to the pair of Riemannian metrics} $%
(h_{11}(t),\varphi _{ij}(x)),$ system which is given by (see the Examples %
\ref{H0} and \ref{G0})%
\begin{equation*}
\frac{d^{2}x^{i}}{dt^{2}}-H_{11}^{1}(t)\frac{dx^{i}}{dt}+\gamma _{jk}^{i}(x)%
\frac{dx^{j}}{dt}\frac{dx^{k}}{dt}=0,
\end{equation*}%
where $H_{11}^{1}(t)$ and $\gamma _{jk}^{i}(x)$ are the Christoffel symbols
of the Riemannian metrics $h_{11}(t)$ and $\varphi _{ij}(x),$ then the
second $h-$KCC-invariant has the form%
\begin{equation*}
\overset{h}{P}\text{ \negthinspace \negthinspace }%
_{j11}^{i}=-R_{pqj}^{i}x_{1}^{p}x_{1}^{q},
\end{equation*}%
where%
\begin{equation*}
R_{pqj}^{i}=\frac{\partial \gamma _{pq}^{i}}{\partial x^{j}}-\frac{\partial
\gamma _{pj}^{i}}{\partial x^{q}}+\gamma _{pq}^{r}\gamma _{rj}^{i}-\gamma
_{pj}^{r}\gamma _{rq}^{i}
\end{equation*}%
are the components of the curvature of the spatial Riemannian metric $%
\varphi _{ij}(x).$ Consequently, the variational equations (\ref%
{Var-Equations}) become the following \textit{jet Jacobi field equations}:%
\begin{equation*}
\frac{\overset{h}{D}}{dt}\left[ \frac{\overset{h}{D}\xi ^{i}}{dt}\right]
+R_{pqm}^{i}x_{1}^{p}x_{1}^{q}\xi ^{m}=0,
\end{equation*}%
where%
\begin{equation*}
\frac{\overset{h}{D}\xi ^{i}}{dt}=\frac{d\xi ^{i}}{dt}+\gamma
_{jm}^{i}x_{1}^{j}\xi ^{m}.
\end{equation*}
\end{example}

\begin{example}
For the particular first order jet rheonomic dynamical system (\ref{Jet_DS})
the jet $h-$deviation curvature $d-$tensor is given by%
\begin{equation*}
\overset{h}{P}\text{ \negthinspace \negthinspace }_{j11}^{i}=\frac{1}{2}%
\frac{\partial ^{2}X_{(1)}^{(i)}}{\partial t\partial x^{j}}+\frac{1}{2}\frac{%
\partial ^{2}X_{(1)}^{(i)}}{\partial x^{j}\partial x^{r}}x_{1}^{r}+\frac{1}{4%
}\frac{\partial X_{(1)}^{(i)}}{\partial x^{r}}\frac{\partial X_{(1)}^{(r)}}{%
\partial x^{j}}+\frac{1}{2}\frac{dH_{11}^{1}}{dt}\delta _{j}^{i}-\frac{1}{4}%
H_{11}^{1}H_{11}^{1}\delta _{j}^{i}.
\end{equation*}
\end{example}

\begin{definition}
The distinguished tensors%
\begin{equation*}
\overset{h}{R}\text{ \negthinspace \negthinspace }_{jk1}^{i}=\frac{1}{3}%
\left[ \frac{\partial \overset{h}{P}\text{ \negthinspace \negthinspace }%
_{j11}^{i}}{\partial x_{1}^{k}}-\frac{\partial \overset{h}{P}\text{
\negthinspace \negthinspace }_{k11}^{i}}{\partial x_{1}^{j}}\right] ,\qquad%
\overset{h}{B}\text{ \negthinspace \negthinspace }_{jkm}^{i}=\frac{\partial 
\overset{h}{R}\text{ \negthinspace \negthinspace }_{jk1}^{i}}{\partial
x_{1}^{m}}
\end{equation*}%
and%
\begin{equation*}
D_{jkm}^{i1}=\frac{\partial ^{3}F_{(1)1}^{(i)}}{\partial x_{1}^{j}\partial
x_{1}^{k}\partial x_{1}^{m}}
\end{equation*}%
are called the \textit{third}, \textit{fourth }and\textit{\ fifth }$h-$%
\textit{KCC-invariant }on the 1-jet vector bundle $J^{1}(\mathbb{R},M)$ of
the system (\ref{SODE}).
\end{definition}

\begin{remark}
Taking into account the transformation rules (\ref{transformations-F}) of
the components $F_{(1)1}^{(i)}$, we immediately deduce that the components $%
D_{jkm}^{i1}$ behave like a $d-$tensor.
\end{remark}

\begin{example}
For the first order jet rheonomic dynamical system (\ref{Jet_DS}) the third,
fourth and fifth $h-$KCC-invariants are zero.
\end{example}

\begin{theorem}[of characterization of the jet $h-$KCC-invariants]
All the five $h-$KCC-invariants of the system (\ref{SODE}) cancel on $J^{1}(%
\mathbb{R},M)$ if and only if there exists a flat symmetric linear
connection $\Gamma _{jk}^{i}(x)$ on $M$ such that%
\begin{equation}
F_{(1)1}^{(i)}=\Gamma _{pq}^{i}(x)x_{1}^{p}x_{1}^{q}-H_{11}^{1}(t)x_{1}^{i}.
\label{h-KCC=0}
\end{equation}
\end{theorem}

\begin{demo}
"$\Leftarrow $" By a direct calculation, we obtain%
\begin{equation*}
\overset{h}{\varepsilon }\text{ }\!\!_{(1)1}^{(i)}=0,\text{ \ \ }\overset{h}{%
P}\text{ \negthinspace \negthinspace }_{j11}^{i}=-\mathfrak{R}%
_{pqj}^{i}x_{1}^{p}x_{1}^{q}=0\text{ and }D_{jkl}^{i1}=0,
\end{equation*}%
where $\mathfrak{R}_{pqj}^{i}=0$ are the components of the curvature of the
flat symmetric linear connection $\Gamma _{jk}^{i}(x)$ on $M.$

"$\Rightarrow $" By integration, the relation%
\begin{equation*}
D_{jkl}^{i1}=\frac{\partial ^{3}F_{(1)1}^{(i)}}{\partial x_{1}^{j}\partial
x_{1}^{k}\partial x_{1}^{l}}=0
\end{equation*}%
subsequently leads to 
\begin{eqnarray*}
\frac{\partial ^{2}F_{(1)1}^{(i)}}{\partial x_{1}^{j}\partial x_{1}^{k}}
&=&2\Gamma _{jk}^{i}(t,x)\Rightarrow \frac{\partial F_{(1)1}^{(i)}}{\partial
x_{1}^{j}}=2\Gamma _{jp}^{i}x_{1}^{p}+\mathcal{U}_{(1)j}^{(i)}(t,x)%
\Rightarrow \\
&\Rightarrow &F_{(1)1}^{(i)}=\Gamma _{pq}^{i}x_{1}^{p}x_{1}^{q}+\mathcal{U}%
_{(1)p}^{(i)}x_{1}^{p}+\mathcal{V}_{(1)1}^{(i)}(t,x),
\end{eqnarray*}%
where the local functions $\Gamma _{jk}^{i}(t,x)$ are symmetrical in the
indices $j$ and $k$.

The equality $\overset{h}{\varepsilon }$ $\!\!_{(1)1}^{(i)}=0$ on $J^{1}(%
\mathbb{R},M)$ leads us to%
\begin{equation*}
\mathcal{V}_{(1)1}^{(i)}=0
\end{equation*}%
and to%
\begin{equation*}
\mathcal{U}_{(1)j}^{(i)}=-H_{11}^{1}\delta _{j}^{i}.
\end{equation*}%
Consequently, we have%
\begin{equation*}
\frac{\partial F_{(1)1}^{(i)}}{\partial x_{1}^{j}}=2\Gamma
_{jp}^{i}x_{1}^{p}-H_{11}^{1}\delta _{j}^{i}
\end{equation*}%
and%
\begin{equation*}
F_{(1)1}^{(i)}=\Gamma _{pq}^{i}x_{1}^{p}x_{1}^{q}-H_{11}^{1}x_{1}^{i}.
\end{equation*}

The condition $\overset{h}{P}$ \negthinspace \negthinspace $_{j11}^{i}=0$ on 
$J^{1}(\mathbb{R},M)$ implies the equalities $\Gamma _{jk}^{i}=\Gamma
_{jk}^{i}(x)$ and%
\begin{equation*}
\mathfrak{R}_{pqj}^{i}+\mathfrak{R}_{qpj}^{i}=0,
\end{equation*}%
where%
\begin{equation*}
\mathfrak{R}_{pqj}^{i}=\frac{\partial \Gamma _{pq}^{i}}{\partial x^{j}}-%
\frac{\partial \Gamma _{pj}^{i}}{\partial x^{q}}+\Gamma _{pq}^{r}\Gamma
_{rj}^{i}-\Gamma _{pj}^{r}\Gamma _{rq}^{i}.
\end{equation*}%
It is important to note that, taking into account the transformation laws (%
\ref{transformations-F}), (\ref{tr-rules-t-s}) and (\ref{rgg}), we deduce
that the local coefficients $\Gamma _{jk}^{i}(x)$ behave like a symmetric
linear connection on $M.$ Consequently, $\mathfrak{R}_{pqj}^{i}$ represent
the curvature of this symmetric linear connection.

On the other hand, the equality $\overset{h}{R}$ \negthinspace \negthinspace 
$_{jk1}^{i}=0$ leads us to $\mathfrak{R}_{qjk}^{i}=0,$ which infers that the
symmetric linear connection $\Gamma _{jk}^{i}(x)$ on $M$ is flat.\medskip
\end{demo}

\textbf{Acknowledgements.} The present research was supported by
the\linebreak Romanian Academy Grant 4/2009.

{\small \setlength{\parskip}{0mm} }

{\small \noindent Authors' addresses: }

{\small \medskip\noindent Vladimir Balan\newline
University Politehnica of Bucharest, Faculty of Applied Sciences,\newline
Department of Mathematics-Informatics I,\newline
Splaiul Independen\c{t}ei 313, RO-060042 Bucharest, Romania.\newline
E-mail: vladimir.balan@upb.ro\newline
Website: http://www.mathem.pub.ro/dept/vbalan.htm }

{\small \medskip \noindent Mircea Neagu\newline
University Transilvania of Bra\c{s}ov, Faculty of Mathematics and
Informatics,\newline
Department of Algebra, Geometry and Differential Equations,\newline
B-dul Iuliu Maniu, Nr. 50, BV 500091, Bra\c{s}ov, Romania.\newline
E-mail: mircea.neagu@unitbv.ro\newline
Website: http://www.2collab.com/user:mirceaneagu }

\end{document}